\documentclass[12pt,twoside]{article}
\usepackage{graphicx, color}

\date{}

\begin{document}

{\emph{Gen. Math. Notes, Vol. 35, No. 1, July 2016, pp.xx-xx}}

\emph{ISSN 2219-7184; Copyright \copyright ICSRS Publication,
2016}

\emph{www.i-csrs.org}

\emph{Available free online at http://www.geman.in}

\centerline{}

\centerline{}

\centerline {\Large{\bf Weak Nonmild Solution of Stochastic}}

\centerline{}

\centerline{\Large{\bf Fractional Porous Medium Equation}}

\centerline{}

%% My definition
\newcommand{\mvec}[1]{\mbox{\bfseries\itshape #1}}

\centerline{\bf {McSylvester Ejighikeme Omaba}}

\centerline{}

\centerline{Department of Mathematics, Computer Science, Statistics and Informatics}
\centerline{Faculty of Science, Federal University Ndufu-Alike Ikwo}
\centerline{PMB 1010 Abakaliki, Ebonyi State, Nigeria}
\centerline{E-mail: mcsylvester\textunderscore omaba@yahoo.co.uk}

\centerline{}

\newtheorem{Theorem}{\quad Theorem}[section]

\newtheorem{Definition}[Theorem]{\quad Definition}

\newtheorem{Corollary}[Theorem]{\quad Corollary}

\newtheorem{Lemma}[Theorem]{\quad Lemma}

\newtheorem{Example}[Theorem]{\quad Example}

\newtheorem{Remark}[Theorem]{\quad Remark}

\newtheorem{Proposition}[Theorem]{\quad Proposition}

\newtheorem{Condition}[Theorem]{\quad Condition}

\newcommand {\N}{\mathrm{N}}
\renewcommand{\P}{\mathrm{P}}
\newcommand {\E}{\mathrm{E}}
\newcommand{\lip}{\mathrm{Lip}}
\newcommand{\Z}{{\mathbf Z}}
\newcommand{\R}{\mathbf{R}}

\centerline{(Received: / Accepted:)}
\centerline{}

\centerline{\bf Abstract}
{\emph{Consider the non-linear stochastic fractional-diffusion equation
 \begin{eqnarray*}
\left \{
\begin{array}{lll}
\frac{\partial}{\partial t}u(x,t)= -(-\Delta)^{\alpha/2}u^m(x,t)+\sigma(u(x,t))\dot{W}(x,t),\, x\in\R^d,t>0,\\
 u(x,0)= u_0(x),\,\,\, x\in\R^d
 \end{array}\right.
\end{eqnarray*}
 with initial data $u_0(x)$ an $L^1(\R^d)$ function, $0<\alpha<2$, and $m>0$. There is no mild solution defined for the above equation because its corresponding heat kernel representation does not exist. We attempt to make sense of the above equation by establishing the existence and uniqueness result via the reproducing kernel Hilbert space (RKHS) of the space-time noise. Our result shows the effect of a space-time white noise on the interaction of fractional operators with porous medium type propagation and consequently studies how the anomalous diffusion parameters influence the energy moment growth behaviour of the system.}}


{\bf 2010 MSC No:} 35R15,\, 35R60,\, 60H15,\,82B44.

%=============================
\section{Introduction}
%=============================
 The authors in \cite{Foondun} studied the following fractional stochastic heat equation
 \begin{eqnarray*}
\left \{
\begin{array}{lll}
 \frac{\partial}{\partial t}u(x,t)= -(-\Delta)^{\alpha/2}u(x,t)+\lambda\sigma(u(x,t))\dot{W}(x,t),\,\,\, x\in\R^d,t>0,\\
 u(x,0)= u_0(x),\,\,\, x\in\R^d.
  \end{array}\right.
\end{eqnarray*}
The operator $-(-\Delta)^{\alpha/2}$ is the fractional Laplacian of order $1 <\alpha\leq 2$. $\lambda$ is a positive parameter called level of noise and $\dot{W}$ denotes space-time white noise. The function $\sigma : \R\rightarrow \R$ is a Lipschitz function satisfying some growth condition. A number of research has been carried out on the above equation for different conditions on the initial function, growth conditions on  $\sigma$ and different types of noises, see \cite{Conus,Foondun,Omaba} and their references.
 We now consider the non-linear fractional-diffusion version of the above equation known as the stochastic fractional porous medium equation (SFPME)
 \begin{eqnarray}\label{eq:FPME}
\left \{
\begin{array}{lll}
 \frac{\partial }{\partial t}u(x,t)= -(-\Delta)^{\alpha/2}u^m(x,t)+\sigma(u(x,t))\dot{W}(x,t),\,\,\, x\in\R^d,t>0,\\
 u(x,0)= u_0(x),\,\,\, x\in\R^d.
  \end{array}\right.
\end{eqnarray}
%\end{document}
  The initial data $u_0(x)$ is an $L^1(\R^d)$ function, the fractional index range $0<\alpha<2$, and the porous medium exponent $m>0$. The difficulty with studying the above non-linear equation is that there is no known heat kernel for the solution and hence one cannot define its mild solution, the reason why little or no extensive work has been done on it; unlike the linear fractional heat equation that have received enough attention, see, \cite{Foondun,Omaba} and their references.  This paper is concerned with the existence, uniqueness and some growth properties/behaviours of the solution using the reproducing kernel Hilbert space of the noise, which is the high point of the research. With the porous medium exponent $m=1$ and $\sigma=0$, is a model of the so-called anomalous diffusion.   See \cite{Vazquez1} which is our main reference paper, for a comprehensive note on the deterministic case of the above equation.  The homogeneous case with $m=1$ has a known integral representation given by $$u(x,t)=\int_{\R^d}p^\alpha(t,x,y)u_0(y)d y,$$ where the kernel $p^\alpha$ has the Fourier transform $\hat{p}^\alpha(t,\xi)=e^{-|\xi|^\alpha t}$. Thus for $0<\alpha<2$, the kernel $p^\alpha$ has the representation $p^\alpha(t,x)=t^{-d/\alpha}F(|x|t^{-1/\alpha})$ for some profile function $F$ that is positive and decreasing and behaves at infinity like $F(r)\sim r^{-(d+\alpha)}$.  That is, the kernel $p^\alpha$ is given explicitly by
$$p^\alpha(t,\,x,\,y) \asymp t^{-d/\alpha}\wedge \frac{t}{|x-y|^{d+\alpha}}\,\,\forall\,\, t>0\,\, x,\,y\in \R^d.$$
There are many physical motivations to considering a fractional Laplacian operator. Modelling diffusive processes using fractional Laplacian has become relevant especially when modelling long-range diffusive interaction. The fractional-diffusion equation of the form $\frac{\partial u}{\partial t}= -(I-\Delta)^{\alpha/2}(-\Delta)^{\gamma/2}(u)+\dot{W}(x,t)$ with $\alpha=2,\,\gamma=0$ arises in neurophysiology. The diffusion operator $(-\Delta)^{1+\gamma},\,\\\gamma>0$, have been used to define hyper-viscosity and to study its effect on the inertial range scaling of fully developed turbulence. Therefore the aim of this article among others is to study the effect of noise on the combination of fractional operators with porous medium type propagation. See \cite{Pablo,Pablo1,Vazquez,Vazquez1}and their references for an in-depth study and analysis of the porous medium equation. The central idea of the paper is to study some growth properties/behaviours of the solution.   This paper is made up of five sections. In section two, we survey some basic distribution concepts that will help to fully understand the paper. Sections three and four give the main results and their proofs respectively, while a brief summary/conclusion of the paper is given in section five.

 \section{Preliminary}
See (\cite{Yosida}, Chapter IX) for more studies on infinitesimal generators and (\cite{Evans}, Chapter 5) for prerequisite concepts and notions used in this paper. The fractional Laplace operator $-(-\Delta)^{\alpha/2},\,\,\alpha\in(0,2)$, are infinitesimal generators of an isotropic stable L\'{e}vy processes. The use of the fractional diffusion operator to replace the normal (standard) Laplace operator is important because it further extends the theory by taking into account the presence of long range interactions.  Let $\Z_{+}=\big\{n\in\Z:n\geq0\big\}$  denote the non-negative integers. A multi index $\alpha=(\alpha_1,...,\alpha_n)\in\Z_+$ is an $n$ tuple of non-negative integers $\alpha_i\geq0$. For multi indices $\alpha=(\alpha_1,...,\alpha_n)$, we define $|\alpha|=\sum_{i=1}^n\alpha_i$. If $x=(x_1,...,x_n)\in\R^n$ and $\alpha=(\alpha_1,...,\alpha_n)\in\Z_+^n$, then we define
 $$D^\alpha=\partial^\alpha=\bigg(\frac{\partial}{\partial x_1}\bigg)^{\alpha_1}...\bigg(\frac{\partial}{\partial x_1}\bigg)^{\alpha_n}=\frac{\partial^{|\alpha|}}{\partial x_1^{\alpha_1}...\partial x_n^{\alpha_1}}$$ as the partial derivative of order $|\alpha|$.
  Let $C^\infty(\Omega)$ be the space of infinitely differentiable function in $\Omega$, that is, the set of functions with continuous partial derivatives of any order and $C_0^\infty(\Omega)$ denote the set of $C^\infty(\Omega)$ functions with compact support in $\Omega$.

 \begin{Definition}(Weak derivatives): Let $u\in L^1(\Omega),\,\,\Omega\subset\R^n$ and $\alpha=(\alpha_1,...,\alpha_n)\in\Z_+^n$. The function $u$ is said to have a weak derivative $D^\alpha u$, if there exists a function $v\in L^1(\Omega)$ such that
 $$\int_\Omega u D^\alpha\varphi d x=(-1)^{|\alpha|}\int_\Omega v\varphi d x,\,\,\varphi\in C_0^\infty(\Omega) $$ where $v=D^\alpha u$.
 \end{Definition}

 \begin{Definition}(Schwartz Space): The Schwartz space, $S(\R^n)$, is the space of functions $\psi \in C^\infty(\R^n)$ such that $$\sup_{x\in\R^n}|x^\alpha\partial_x^\beta\psi(x)|<\infty,\,\,\,\forall\,\,\alpha,\,\beta\in\Z_+^n.$$ Suppose $\varphi\in S,$ then for every $d\in\N$ and $\alpha\in\Z_+^n$, there exists a constant $C_{d,\alpha}$ such that $$|\partial^\alpha\varphi(x)|\leq\frac{C_{d,\alpha}}{(1+|x|^2)^{\frac{d}{2}}},\,\,\forall\,x\in\R^n.$$
 \end{Definition}

\begin{Definition}The space of tempered distribution, $S'(\R^d)$, is the space of linear continuous functional on $S(\R^d)$. That is, $T\in S'(\R^d)$ if the following conditions hold:
 \begin{itemize}
 \item Linearity: $T(\alpha\psi_1+\beta\psi_2)=\alpha T(\psi_1)+\beta T(\psi_2)$
 \item Continuity: $T(\psi_n)\rightarrow T(\phi)$, if $\phi_n\rightarrow\phi$ as $n\rightarrow\infty,\,\,\psi_n,\psi\in S(\R^d).$
 \end{itemize}
 \end{Definition}
 
 Assume that $\gamma$ is a non-negative integer. Then we define the Sobolev space $H^\gamma(\R^d)$, to be the set of all $\psi\in L^2(\R^d)$ whose (distributional) derivatives $\partial^\alpha\psi$ belong to $ L^2(\R^d)$ for $|\alpha|\leq\gamma$:
 $$H^\gamma(\R^d)=\bigg\{\psi\in L^2(\R^d)\,\,\,|\,\,\,\partial^\alpha\psi\in L^2(\R^d)\,\, \textrm{for}\,\, |\alpha|\leq\gamma\bigg\}.$$ For fractional and negative exponents, the integral operators and their inverses can be defined as bounded operators on the fractional Sobolev spaces. Thus, we give an equivalent characterisation of this definition in terms of Fourier transforms as follows.
 
 \begin{Definition} For $\gamma\in\R^d$, the Sobolev space $H^\gamma(\R^d)$ is defined by
 $$H^\gamma(\R^d)=\bigg\{\psi\in S'(\R^d)\,\,\,|\,\,\,(1+|\xi|^2)^{\gamma/2}\hat{\psi}\in L^2(\R^d)\bigg\},$$ where $(1+|\xi|^2)^{\alpha/2}$ is the characteristic polynomial of $(I-\Delta)^{\gamma/2}$.
 \end{Definition}
 
 The inner product of two functions $\psi$ and $\varphi$ in $H^\gamma(\R^d)$ is defined as
 $$\langle\psi,\varphi\rangle_{H^\gamma(\R^d)}=\langle (1+|\xi|^2)^{\gamma/2}\hat{\psi},\, (1+|\xi|^2)^{\gamma/2}\hat{\varphi}\rangle_{L^2(\R^d)},$$ and its norm $$\|\psi\|_{H^\gamma(\R^d)}=\|(1+|\xi|^2)^{\gamma/2}\hat{\psi}\|_{L^2(\R^d)}.$$
   The fractional Laplace operator $(-\Delta)^{\alpha/2}$, is defined through Fourier transform as follows: let $\psi$ be a test function and $(-\Delta)^{\alpha/2}\psi=\varphi$, then
 \begin{equation}\label{eq:Fourier}
 \hat{\varphi}(\xi)=|\xi|^\alpha\hat{\psi}(\xi).
 \end{equation}
   Gaussian space-time white noise belongs to a special class of Sobolev space known as a local Sobolev space, see \cite{Walsh}.

\begin{Definition} A function $u(x,t)$ is locally square integrable on $\R^d\times (0,\infty)$ if
$$\int_D |u(x,t)|^2 d x d t<\infty,$$ for every compact set $D$ in $\R^d\times (0,\infty)$.
 \end{Definition}
 
 \begin{Definition} Suppose $\gamma\in\R^d$. If $u$ is any distribution, we say that $u\in H^\gamma_{loc}(\R^d\times (0,\infty))$, if for any $\psi(x,t)\in C_0^\infty(\R^d\times (0,\infty))$, \,\,\,$u(x,t)\psi(x,t)\in  H^\gamma(\R^d\times (0,\infty)).$
 \end{Definition}
 There is no known kernel for the equation (\ref{eq:FPME}) and therefore we cannot define a mild solution to the equation, so we give the weak solution in terms of test functions.  Given that $\varphi$ and $\psi$ are test functions, then it follows by equation (\ref{eq:Fourier}) with Plancherel's identity that
  \begin{eqnarray*}\int_{\R^d}(-\Delta)^{\alpha/2}\varphi\psi d x=\int_{\R^d}|\xi|^\alpha\hat{\varphi}\hat{\psi} d\xi&=&\int_{\R^d}|\xi|^{\alpha/2}\hat{\varphi}|\xi|^{\alpha/2}\hat{\psi} d\xi\\&=&\int_{\R^d}(-\Delta)^{\alpha/4}\varphi(-\Delta)^{\alpha/4}\psi d x.\end{eqnarray*}
 We begin by first considering when $\sigma=1$:
 \begin{eqnarray}\label{eq:FPME:2}
\left \{
\begin{array}{lll}
 \frac{\partial}{\partial t}u(x,t)= -(-\Delta)^{\alpha/2}u^m(x,t)+\dot{W}(x,t),\qquad x\in\R^d,t>0,\\
 u(x,0)= u_0(x),\,\,\, x\in\R^d.
  \end{array}\right.
\end{eqnarray}
Before making sense of the above equation, we state the following 
proposition:

\begin{Proposition}[\cite{Walsh}, Proposition 9.5]
Let $\epsilon>0$ and $\dot{W}(x,t)$ a Gaussian space-time white noise on $\R^d\times[0,\infty)$. Then with probability one, $$\dot{W}(x,t)\in  H^{-d/2-\epsilon}_{loc}(\R^d\times(0,\infty)).$$
\end{Proposition}

Now, multiply equation (\ref{eq:FPME:2}) by a test function $\psi$ and integrate by part,
$$\int_0^T\int_{\R^d}u_t\varphi d x d t=-\int_0^T\int_{\R^d}(-\Delta)^{\alpha/2}u^m\varphi d x d t+\int_0^T\int_{\R^d}\dot{W}(x,t)\varphi d x d t.$$
$$-\int_0^T\int_{\R^d}u\varphi_t d x d t=-\int_0^T\int_{\R^d}(-\Delta)^{\alpha/4}u^m(-\Delta)^{\alpha/4}\varphi d x d t-\int_0^T\int_{\R^d}W(x,t)\varphi_t d x d t.$$ In particular, we consider $W$ a Wiener process in $S'(\R^d)$ whose spectral measure is equal to $\delta_0$. The process $W$ has the form $W(x,t)=(2\pi)^{-d/4}B(t),\,\,\,x\in\R^d,\,\,t\geq0$, where $B$ is a real-valued Wiener process, see (\cite{Peszat}, Example 14.21). Then \begin{eqnarray}\label{eq:main}
\int_0^T\int_{\R^d}u(x,t)\varphi_t(x,t) d x d t=\int_0^T\int_{\R^d}(-\Delta)^{\alpha/4}u^m(x,t)(-\Delta)^{\alpha/4}\varphi d x d t\hskip0.0000000000in\\+(2\pi)^{-d/4}\int_0^T\int_{\R^d}B(t)\varphi_t(x,t) d x d t.\nonumber
\end{eqnarray}

\section{Main Results}

On the space $S'(\R^d)$, $W$ is a square integrable mean-zero L\'{e}vy process that can be regarded as a square integrable L\'{e}vy process on a properly chosen Hilbert space. So we give the RKHS of the process $W$ in terms of its covariance $K$. Define $\mathcal{H}$ as the set of all $W\in S'(\R^d)$ such that
\begin{equation}\label{eq:spatial}|\langle W,\varphi\rangle|\leq L\sqrt{K(\varphi,\varphi)},\,\,\,\,\forall\,\varphi\in S(\R^d),\end{equation} with a constant $L<\infty$ independent of $\varphi$. Let $(H_n,\langle.,.\rangle_{H_n}),\,\,n\in\N$, be a decreasing sequence of separable Hilbert spaces, then we have:

\begin{Lemma}[\cite{Peszat}, Lemma 14.9]\label{eq:inequality} There exist $n$ and $C$ such that $|K(\varphi,\psi)|\leq C|\varphi|_{H_n}|\psi|_{H_n}$ for all $\varphi,\psi\in S(\R^d)$.
\end{Lemma}

The next theorem states that $W$ can be seen as a spatially stationary random field on $\R^d\times[0,\infty)$:

\begin{Theorem}[\cite{Peszat}, Theorem 14.24]\label{eq:Peszat}
If the spectral measure $\mu$ of $W$ is finite, then $W$ can be identified with a random field $W(x,t),\,\,t\geq0,\,\,x\in\R^d,$ that is, $$\langle W(t),\psi(t)\rangle=\int_{\R^d}W(x,t)\psi(x,t) d x,\,\,t\geq0,\,\,\psi\in S(\R^d\times[0,\infty)).$$
\end{Theorem}

We now return to equation (\ref{eq:FPME:2}). We multiply equation (\ref{eq:FPME:2}) by a test function $\psi$ and integrate by part, thus

\begin{eqnarray}\label{eq:main}
\int_0^T\int_{\R^d}u(x,t)\psi_t(x,t)d x d t=\int_0^T\int_{\R^d}(-\Delta)^{\alpha/4}u^m(x,t)(-\Delta)^{\alpha/4}\psi d x d t\hskip0.0000000000in\\+\int_0^T \langle W(t),\psi_t(t)\rangle d t.\nonumber
\end{eqnarray}
We can make sense of the above integrals if $u$ and $u^m$ are to lie on a well defined spaces, and the possible suitable space for $u^m$ is the fractional Sobolev space $H^{\alpha/2}(\R^d)$, defined as the completion of $C_0^\infty(\R^d)$ with the following norm:
$$\|\varphi\|_{H^{\alpha/2}(\R^d)}=\bigg(\int_{\R^d}|\xi|^\alpha|\hat{\varphi}|^2 d\xi\bigg)^{1/2}=\|(-\Delta)^{\alpha/4}\varphi\|_{L^2(\R^d)}.$$

\begin{Definition}A function $u$ is a weak solution to equation (\ref{eq:FPME:2}) if:
\begin{enumerate}
\item[1.] $u\in L^1(\R^d\times(0,T))$ for all $T>0$,\,$u^m\in L^2_{loc}(H^{\alpha/2}(\R^d);(0,\infty))$,
\item[2.] $W$ satisfies equation (\ref{eq:spatial}) with Lemma \ref{eq:inequality}.
\item[3.]equation (\ref{eq:main}) holds for every $\varphi\in C_0^1(\R^d\times (0,T))$,
\item[4.] $u(.,t)\in L^1(\R^d)$ for all $t>0$, $\lim_{t\rightarrow 0^+}u(.,t)=u_0\in L^1(\R^d)$.
\end{enumerate}
\end{Definition}

Following the above known results and definition, we thus give main results of this article.

\begin{Proposition}\label{eq:property}Suppose $u_0\in L^1(\R^d)$, and let $u\in L^1(\R^d\times[0,T])$ for all $T>0$ be a weak solution to (\ref{eq:FPME}),\,$u^m\in L^2_{loc}(H^{\alpha/2}(\R^d);(0,\infty))$, and $W$ satisfies equation (\ref{eq:spatial}) with Lemma \ref{eq:inequality} then $\int_{\R^d}u(x,t)d x\in C(\R^+)$, that is, $u\in C(L^1(\R^d);[0,T])$ and
\begin{eqnarray*}\int_{\R^d}u(x,t)d x&=&\int_0^t\int_{\R^d}\sigma(u(x,s))\dot{W}(x,s)d x d s+\int_{\R^d}u_0(x)d x\\&=&\int_0^t\int_{\R^d}\sigma(u(x,s))W(d x,d s)+\int_{\R^d}u_0(x)d x.\end{eqnarray*}
\end{Proposition}

We will give the proof of proposition \ref{eq:property} for $\sigma=1$ for simplicity, but the same steps apply when $\sigma\neq1$, with $\sigma:\R\rightarrow\R$ Lipschitz continuous. For existence and uniqueness, we need the following condition on $\sigma$. Essentially this condition says that $\sigma$ is globally Lipschitz in the first variable.

\begin{Condition}\label{cond:E-U}
There exist a finite positive constant, $\lip_\sigma$ such that for all $x,\,y\in \R$, we have
$\sigma(0)=0,$\,\, and $|\sigma(x)-\sigma(y)|\leq \lip_\sigma|x-y|.$
\end{Condition}

The existence and uniqueness result of the main equation follows by equation (\ref{eq:spatial}), Lemma \ref{eq:inequality} and Theorem \ref{eq:Peszat}.

\begin{Theorem}\label{eq:ex-un}
Under condition \ref{cond:E-U}, there exists a unique solution to (\ref{eq:FPME:2}).
\end{Theorem}

The next result is the growth property of the solution using proposition \ref{eq:property}.

\begin{Theorem}[Exponential growth]\label{exponential} There exist constants $\lip_\sigma>0$ and $C>0$ such that $$\|u(t)\|_{L^2(\P)}\leq\frac{1}{\sqrt{C}}\|u_0\|_{L^1(\R^d)}\exp\bigg(\frac{\lip_\sigma}{\sqrt{C}} t\bigg).$$
\end{Theorem}

\section{Proofs of Main Results}We now give proofs to our main results:\\

{\bf Proof of Proposition \ref{eq:property}:}
\\\\
Let's define a test function $\psi$ by $\psi(x,t)=\varphi_R(x)\phi(t)$ where $\varphi_R(x)=\varphi(\frac{x}{R})$ with the cut-function $\varphi$ defined as follows:
$\varphi\in C_0^\infty(\R^d),\,\,\,\, 0\leq \varphi(.)\leq 1,\,\,\varphi(x)=1$\,\,\,for\,\,$|x|\leq 1,\,\,\,\,\varphi(x)=0$\,\,\,for $\,\,\,|x|\geq 2.$
We have for any test function $\phi\in C_0^\infty(\R^+)$,
\begin{eqnarray*}-\int_0^\infty\int_{\R^d}u(t,x)\varphi_R(x)\phi'(t) d t d x\hskip7.9in\\=\int_0^\infty\int_{\R^d}\bigg[ -u^m(x,t)(-\Delta)^{\alpha/2}\varphi_R(x)+\dot{W}(x,t)\varphi_R(x)\bigg]\phi(t)d t d x\hskip5.2in\\+\int_{\R^d}u_0(x)\varphi_R(x)\phi(0)d x\hskip5.25in\\
=\int_0^\infty\int_{\R^d}\bigg[- u^m(x,t)\frac{1}{R^\alpha}(-\Delta)^{\alpha/2}\varphi(\frac{x}{R})+\dot{W}(x,t)\varphi_R(x)\bigg]\phi(t)d t d x\hskip5.0in\\ +\int_{\R^d}u_0(x)\varphi_R(x)\phi(0)d x.\hskip5.0in
\end{eqnarray*}
Taking limit as $R\rightarrow\infty$, we obtain using the Lebesgue dominated convergence theorem, that $$-\int_0^\infty\int_{\R^d}u(x,t)\phi'(t)d t d x=\int_0^\infty\int_{\R^d}\dot{W}(x,t)\phi(t)d t d x+\int_{\R^d}u_0(x)\phi(0)d x.$$
 Let $F_u(t)=\int_{\R^d}u(x,t)d x$ and $F_{\dot{W}}(t)=\int_{\R^d}\dot{W}(x,t)d x$ where $F_u(t)$ belong to $L_{loc}^1(\R^+)$ and $F_{\dot{W}}(t)$ belong to $L_{loc}^2(\P)$. Let $T>0$ and consider the test function $\phi(t)=1_{\{0\leq t\leq T\}}$; also choose a sequence of test functions $\phi_n\in C_0^\infty(\R^+)$ with $\phi_n$ decreasing, $\phi_n(t)\leq\phi(t)$ and $\phi_n(t)=1$ on the interval $[0,T-1/n]$ for sufficiently large $n$. Then we have the following by applying Lebesgue dominated convergence theorem on $F_\sigma$:
 \begin{eqnarray*}-\int_0^\infty F_u(t)\phi_n'(t)d t\hskip7.7in\\=\int_0^\infty F_{\dot{W}}(t)\phi_n(t)d t+\int_{\R^d}u_0(x)\phi_n(0)d x\rightarrow\int_0^T F_{\dot{W}}(t)d t+\int_{\R^d}u_0(x)d x\hskip4.1in\end{eqnarray*} as $n\rightarrow\infty$ and $\delta\rightarrow 0$. This shows that $F_u\in C(\R^+)$ or better still, has a continuous representation in its Lebesgue class. Now choosing $\phi_n$ nicely, at each Lebesgue point $T$ of $F_u$, one also passes a limit on the left hand side to obtain that $$F_u(T)=\int_0^T F_{\dot{W}}(t)d t+\int_{\R^d}u_0(x)d x.$$
\\\newpage
{\bf Proof of Theorem \ref{eq:ex-un}:}
\\\\
By Lemma \ref{eq:inequality}, equation (\ref{eq:spatial}) and H\"{o}lder's inequality, we have that
\begin{eqnarray*}
\bigg\arrowvert\int_0^T\int_{\R^d}u(x,t)\psi_t(x,t)d x d t\bigg\arrowvert \hskip3.0000000000in\\\leq\|(-\Delta)^{\alpha/4}u^m\|_{L^2(\R^d\times(0,T))}.\|(-\Delta)^{\alpha/4}\psi\|_{L^2(\R^d\times(0,T))}+ T\,L\sqrt{C}\sup_{0\leq t\leq T}|\psi_t|_{H_n}\\
=\|u^m\|_{H^{\alpha/2}(\R^d\times(0,T))}.\|\psi\|_{H^{\alpha/2}(\R^d\times(0,T))}+ T\,L\sqrt{C}|\psi|_{H_n}<\infty.\hskip0.910000000000in
\end{eqnarray*}
Going back to equation (\ref{eq:FPME}), we have that
\begin{eqnarray*}
\E\bigg\arrowvert\int_{\R^d}[u_1(x,t)-u_2(x,t)]d x\bigg\arrowvert^2&=&\E\bigg\arrowvert\int_0^t\int_{\R^d}[\sigma(u_1(x,t))-\sigma(u_2(x,t))]W(d x, d t)\bigg\arrowvert^2\\&\leq&\lip_\sigma^2\E\int_0^t\int_{\R^d}\big\arrowvert u_1(x,t)-u_2(x,t)\big\arrowvert^2 d x d t.
\end{eqnarray*} Let there exists $C>0$ such that $$C\int_{\R^d}\E|u_1(x,t)-u_2(x,t)|^2d x\leq\E\bigg\arrowvert\int_{\R^d}[u_1(x,t)-u_2(x,t)]d x\bigg\arrowvert^2,$$ then it follows that $\|u_1-u_2\|_{L^2(\P)}\leq \lip_\sigma\sqrt{T/C}\|u_1-u_2\|_{L^2(\P)}$ which proves the uniqueness given that $1-\lip_\sigma\sqrt{T/C}>0$.
\\\\
{\bf Proof of Theorem \ref{exponential}:}
\\\\
From the above proposition \ref{eq:property}, we have that
\begin{eqnarray*}
\|u(t)\|_{L^2(\P)}\leq\frac{\lip_\sigma}{\sqrt{C}} \int_0^t\|u(s)\|_{L^2(\P)}d s+\frac{1}{\sqrt{C}}\|u_0\|_{L^1(\R^d)},
\end{eqnarray*}and the result follows by Gronwall's lemma.
%
%==========================
\section{Conclusion}
%==========================
We are able to establish the not so easy existence and uniqueness of the solution by reproducing kernel Hilbert space of the noise since there does not exist a heat kernel and consequently no mild solution to the non-linear equation. Also proved is the second moment growth estimate for the solution.

\end{document}